\documentclass[a4paper]{amsart}
\usepackage{amsmath,amssymb,amsthm}
\usepackage{hyperref}

\newcommand{\g}{\mathfrak{g}}
\newcommand{\h}{\mathfrak{h}}
\newcommand{\m}{\mathfrak{m}}

\newcommand{\so}{\mathfrak{so}}

\newcommand{\gl}{\mathfrak{gl}}

\newcommand{\stab}{\mathfrak{stab}}

\newcommand{\R}{\mathbb{R}} 
\newcommand{\RP}{\mathbb{R}P} 

\newcommand{\dcf}{\mathcal{C}} 
\newcommand{\bggCup}{\mathbin{\scriptstyle\sqcup}} 

\newcommand{\cD}{\mathcal{D}}
\newcommand{\cN}{\mathcal{N}}

\newcommand{\tR}{\underline{\R}}

\newcommand{\del}{\partial}
\newcommand{\quabla}{\square}

\DeclareMathOperator{\im}{im}
\DeclareMathOperator{\ad}{ad}

\DeclareMathOperator{\id}{id}
\DeclareMathOperator{\Hom}{Hom}
\DeclareMathOperator{\GL}{GL}

\newcommand{\HkW}{H_k(T^*\Sigma,W)} 
\newcommand{\cs}{\mathsf{c}} 
\newcommand{\hR}{\hat{R}} 
\newcommand{\trace}{\mathrm{tr}}
\newcommand{\dDg}{d^{\cD^g}} 
\newcommand{\cDg}{\cD^g} 
\newcommand{\cNg}{\cN^g} 
\newcommand{\RDg}{R^{\cD^g}} 
\newcommand{\dgt}{d^g_t} 
\newcommand{\CYc}{C} 

\newcommand{\Cech}{\v{C}ech } 

\newcommand{\tg}{\underline{\g}}

\newcommand{\trR}{\underline{\R}}

\newcommand{\Cinfty}{\mathrm{C}^\infty}
\newcommand{\fsec}[1]{\Cinfty(\submfd,#1)} 
\newcommand{\submfd}{\Sigma}
\newcommand{\fone}[1]{\Omega^1(\submfd,#1)} 
\newcommand{\fk}[1]{\Omega^k(\submfd,#1)} 
\newcommand{\vspan}[1]{\langle#1\rangle}
\newcommand{\clq}{\tilde{g}}

\theoremstyle{plain}
\newtheorem{theorem}{Theorem}[section]
\newtheorem{proposition}[theorem]{Proposition}

\theoremstyle{remark}
\newtheorem{remark}[theorem]{Remark}

\theoremstyle{definition}
\newtheorem{definition}[theorem]{Definition}

\begin{document}

\title{M\"obius-flat hypersurfaces in projective space}
\author{Daniel J. Clarke}
\email{daniel.jc@gmail.com}

\begin{abstract}
I give a theory of M\"obius-flat hypersurfaces in $\RP^n$, analogous to that in conformal $S^n$. This unifies the classes of hypersurfaces with flat induced conformal structure ($n > 3$) and a classically studied class of surfaces ($n = 3$). I extend an example of Akivis--Konnov, and use polynomial conserved quantities to characterise hypersurfaces with flat centro-affine metric among M\"obius-flat hypersurfaces. Part of the theory has an obvious counterpart in Lie sphere geometry.
\end{abstract}

\maketitle

\section{Introduction} \label{sec:intro}
This paper concerns two classes of projective hypersurface. Around the late 1920's, \Cech \cite{vCech1928, vCech1929} and Kauck\'y \cite{Kauck'y1932} studied a class of hypersurfaces in $\RP^n$ possessing a one-parameter family of asymptotic deformations rescaling the Darboux cubic form. 
Slightly more recently, the conformal structure induced on a projective hypersurface with non-degenerate second fundamental form (``tangentially non-degenerate'' in the language of Akivis--Konnov \cite{Akivis1993}) was studied, and in particular those that are conformally flat.

In the realm of conformal geometry, Burstall--Calderbank \cite{Burstall2010,Burstallb} unify conformally flat submanifolds having flat normal bundle in conformal $S^n$, with Guichard and channel surfaces by taking an integrable systems viewpoint. Furthermore, they provide a \emph{conformal} approach to constant Gaussian curvature submanifolds of spaceforms. This is achieved using the concept of a \emph{polynomial conserved quantity}. These have also been used by Burstall--Santos \cite{Burstall2012} to view constant mean curvature hypersurfaces in spaceforms as isothermic surfaces, Quintino \cite{Quintino2008} in the context of Willmore surfaces, and Burstall--Calderbank in their study of Guichard surfaces as Lie applicable surfaces (private communication).

 Classically, surfaces in projective space were studied \cite{Wilczynski1907} via solutions of the linear system
\begin{align*}\\
\sigma_{xx} &= \beta\sigma_y + \frac{1}{2}(V - \beta_y)\sigma\\
\sigma_{yy} &= \gamma\sigma_x + \frac{1}{2}(W - \gamma_x)\sigma.
\end{align*}
where $x,y$ are asymptotic co-ordinates (we adopt here the notation of \cite{Ferapontov2004}). Kauck\'y \cite{Kauck'y1932} and \Cech \cite{vCech1928,vCech1929} identified surfaces satisfying
\begin{equation}
	\label{eq:cmf}
	\beta_{yyy} = \gamma_{xxx}
\end{equation}
as possessing a one-parameter family of asymptotic deformations rescaling the Darboux cubic form $\dcf$.

The condition \eqref{eq:cmf} depends upon the choice of asymptotic co-ordinates, so instead we prefer to work with the condition that there exists functions $a,b$ so that
\begin{equation}
	\label{eq:cmf_a}
	2\beta_y b - \beta b_y = 2\gamma_x a - \gamma a_x
\end{equation}
and
\begin{align}
	2b_x = 2\gamma\beta_y + \beta\gamma_y
	\label{eq:cmf_b}\\
	2a_y = 2\beta\gamma_x + \gamma\beta_x
	\label{eq:cmf_c}
\end{align}
 
Upon change of co-ordinates, the functions $a, b$ scale as components of a quadratic differential. We shall show how equations \eqref{eq:cmf_a}-\eqref{eq:cmf_c} may be expressed in terms of this quadratic differential using a differential pairing of Calderbank--Diemer \cite{Calderbank2001} in \S\ref{sec:hom_surfaces}.

 We will introduce a notion of projectively M\"obius-flat hypersurface in terms of the existence of one-forms satisfying certain algebraic and differential conditions. In particular, we obtain a zero-curvature formulation for conformal flatness of projective hypersurfaces.

After studying this notion in its own right, we will show that it unifies hypersurfaces with flat asymptotic conformal structure when $n > 3$, and surfaces satisfying  \eqref{eq:cmf_a}-\eqref{eq:cmf_c} when $n=3$.

In centro-affine geometry, i.e., differential geometry in a finite-dimensional vector space, Ferapontov \cite{Ferapontov2004} studied hypersurfaces with flat centro-affine metric.  He showed that these hypersurfaces possess a spectral deformation. Also he found that hypersurfaces with flat centro-affine metric satisfy \eqref{eq:cmf} when $n=3$, observing that they are conformally flat (in the projective sense) when $n > 3$. Otherwise said, hypersurfaces with flat centro-affine metric are projectively M\"obius-flat.

In this paper we will see how this fits into the framework of polynomial conserved quantities. For surfaces, the existence of a polynomial conserved quantity may be formulated as
\begin{theorem}
	\label{thm:pcq_classical}
Let $\alpha$ be a potential for the Chebyshev covector, and suppose $\beta_{xxx} = \gamma_{yyy}$. Then the surface has flat centro-affine metric if and only if
\begin{align}
\beta_y&=2\alpha_{xx} \\
\gamma_x&=2\alpha_{yy}\\
V &= 2(\beta\alpha_y + \alpha_x^2)\\
W &= 2(\gamma\alpha_x + \alpha_y^2)\\
1 &= \beta\gamma -  4\alpha_x\alpha_y
\end{align}
Moreover, if $\alpha$ is any function satisfying these equations then $\alpha$ is a potential for the Chebyshev covector of a projective transform of a surface with flat centro-affine metric.
\end{theorem}
but we will derive this from a more natural formulation that holds for $m\geq 3$ also.

Thus, the main objectives of this paper are: to provide an account of M\"obius-flat projective hypersurfaces analogous to that in conformal geometry, and to characterise among these the hypersurfaces with flat centro-affine metric.
It is evident that the work in \S\S\ref{sec:definition},\ref{sec:conformally_flat} has a direct analogue for surfaces in Lie sphere geometry using Lie's line-sphere correspondence.

In \S\ref{sec:newbit} we will review the ``gauge-theoretic'' formalism that we use for projective and conformal differential geometry. We will then examine the known spectral deformation to provide motivation for our definition of projectively M\"obius flt hypersurface which we introduce in \S\ref{sec:definition}. This is followed by two examples in \S\ref{sec:examples}: second-order envelopes of quadric congruences (dealt with already by Akivis-Konnov for $m>4$) and hypersurfaces with flat centro-affine metric (dealt with by Ferapontov as mentioned above). We have postponed to \S\ref{sec:conformally_flat} the matter of proving that our definition in \S\ref{sec:definition} unifies the two classes as claimed above. The remainder of this section and \S\ref{sec:pcq_surface} then covers the special case of surfaces, in an attempt to make all this more concrete for readers who are used to working in the setup of equations \eqref{eq:cmf_a}-\eqref{eq:cmf_c}.

This work is based on part of the author's PhD thesis.
\section{Notation and motivation} \label{sec:newbit}

First we recall the gauge-theoretic setup for projective differential geometry from \cite{Burstall2004}.

Let $\submfd$ be a fixed $m$-manifold; in this paper we always take $m\geq 2$. For a vector space $V$, we use $\underline{V}$ to denote the trivial vector bundle $V\times\submfd$ over $\submfd$. This bundle carries a trivial connection $d$ given by $d_X\sigma:=(d\sigma)(X)$ where the second $d$ is just the usual exterior derivative. Given a vector bundle $E$, we denote by $\fsec{E}$ and $\fk{E}$ the spaces of sections of $E$ and $k$-forms with values in $E$ respectively. 

Throughout we will identify subbundles of $\trR^{n+1} := \submfd\times\R^{n+1}$ with maps from $\submfd$ into Grassmannians. More precisely, a rank $k$ subbundle $U\subseteq\trR^{n+1}$ corresponds to the (smooth) map $x\mapsto U_x$ from $\submfd$ to the Grassmannian of $k$-dimensional subspaces of $\R^{n+1}$, where $U_x$ is the fibre of $U$ over the point $x$. In particular let $\Lambda$ be a codimension one immersion of $\submfd$ into $\RP^n$. The \emph{contact lift} is then a filtration $E$,
\begin{equation*}
E_0 = 0 \subseteq E_1 = \Lambda \subseteq E_2 = d\Lambda \subseteq E_3 = \trR^{n+1}
\end{equation*}
where $d\Lambda$ is a shorthand we use to denote the subbundle whose sections are of the form $d_X\sigma$ for $X\in\fsec{TM}$ and $\sigma\in\fsec{\Lambda}$.

In this paper we will take the view that a surface is a rank $(1,n)$ filtration $E$ by subbundles of $\trR^{n+1}$ satisfying $dE_1 = E_2$; such $E$ we call \emph{Legendre}. This induces a filtration of the bundle of Lie algebras $\tg := \submfd\times\gl(\R^{n+1})$ with
\begin{align*}
\tg_{-2} &= \{X \in \tg | XE_2 = 0, X\trR^{n+1}\subseteq E_1\} \\
\tg_{-1} &= \{X \in \tg | XE_1 = 0, XE_2\subseteq E_1\} \\
\tg_{0} &= \{X \in \tg | XE_1\subseteq E_1, XE_2\subseteq E_2\} \\
\tg_{1} &= \{X \in \tg | XE_1\subseteq E_2\}
\end{align*}

Given a metric $g$ on $\trR^{n+1}$, there is a unique decomposition $d = \cDg + \cNg$ where $\cDg$ is a connection satisfying $\cDg g = 0$ and $\cNg$ is a $\tg$-valued one-form whose values are $g$-symmetric. These are the $\h$ and $\m$ parts respectively of $d$, where $\h:=\{X\in\tg|\text{$X$ skew for g}\}$ and $\m:=\{X\in\tg|\text{$X$ symmetric for g}\}$. The following formulae are useful for calculation:
\begin{align*}
g(\cD^g\sigma,\tau) &= \frac{1}{2}(g(d\sigma,\tau) - g(\sigma,d\tau) + dg(\sigma,\tau))\\
g(\cN^g\sigma,\tau) &= \frac{1}{2}(g(d\sigma,\tau) + g(\sigma,d\tau) - dg(\sigma,\tau)),\\
\end{align*}
for $\sigma,\tau\in\fsec{\trR^{n+1}}$. We say $E$ \emph{envelopes} $g$ if $E_1$ is null for $g$ and $\cNg\in\fone{\tg_{0}}$; we say $g$ is \emph{unimodular} if the values of $\cNg$ are trace-free. We will recall more on this in \S\ref{sec:examples}.

Let $q$ be a quadric congruence in $\RP^n$. This is the same thing as a rank one subbundle of $S^2(\trR^{n+1})^*$ all of whose sections are non-degenerate. There is then a distinguished metric in this conformal class (up to constant rescaling) satisfying $\trace\cN^g=0$, which we call \emph{unimodular}.

Now let us see how this can help us analyse the spectral deformation associated with \eqref{eq:cmf}.

We calculate
\begin{align}
\sigma_{xxy} &= \beta_y\sigma_y + \beta(\gamma\sigma_x + \frac{1}{2}(W-\gamma_x)\sigma) + \frac{1}{2}(V_y - \beta_{yy})\sigma + \frac{1}{2}(V - \beta_y)\sigma_y\\
\sigma_{yyx} &= \gamma_x\sigma_x + \gamma(\beta\sigma_y + \frac{1}{2}(V-\beta_y)\sigma) + \frac{1}{2}(W_x - \gamma_{xx})\sigma + \frac{1}{2}(W - \gamma_x)\sigma_x.
\end{align}

Thus,
\begin{align}
\begin{pmatrix}
\sigma\\
\sigma_x\\
\sigma_y\\
\sigma_{xy}
\end{pmatrix}_x
&=
\begin{pmatrix}
0 & 1 & 0 & 0 \\
\frac{1}{2}(V - \beta_y) & 0 & \beta & 0 \\
0 & 0 & 0 & 1 \\
\frac{1}{2}(\beta W - \beta\gamma_x + V_y - \beta_{yy}) & \beta\gamma & \frac{1}{2}(V + \beta_y) & 0 
\end{pmatrix}
\begin{pmatrix}
\sigma\\
\sigma_x\\
\sigma_y\\
\sigma_{xy}
\end{pmatrix}\\
\begin{pmatrix}
\sigma\\
\sigma_x\\
\sigma_y\\
\sigma_{xy}
\end{pmatrix}_y
&=
\begin{pmatrix}
0 & 0 & 1 & 0 \\
0 & 0 & 0 & 1 \\
\frac{1}{2}(W - \gamma_x) & \gamma & 0 & 0 \\
\frac{1}{2}(\gamma V - \gamma\beta_y + W_x - \gamma_{xx})  & \frac{1}{2}(W + \gamma_x) & \beta\gamma & 0 
\end{pmatrix}
\begin{pmatrix}
\sigma\\
\sigma_x\\
\sigma_y\\
\sigma_{xy}
\end{pmatrix}.
\end{align}
Putting $\psi = (\sigma, \sigma_x, \sigma_y, \sigma_{xy})$, this is of the form $d\psi = A\psi$ where $A$ is a matrix-valued one-form. The one-form $A$ represents the connection $d$ w.r.t. the frame $\psi$. 

We now insert a spectral parameter into $A$ by multiplying $\beta, \gamma$ by $t$, and $V,W$ by $t^2$ (see e.g. \cite{Ferapontov2004}), and call the result $A_t$. This defines a new connection according to $d_t\psi = A_t\psi$. The difference $d_t - d$ between these connections is a Lie algebra-valued one-form, which in matrix form w.r.t. the frame $\psi$ is $(A_t - A_1)^T$. 

Let us write this out explicitly, separating the powers of $t$:
\begin{align*}
\frac{\partial}{\partial x}(A_t - A_1)^T &= (t - 1)
\begin{pmatrix}
0 & -\frac{1}{2}\beta_y & 0 & -\frac{1}{2}\beta_{yy} \\
0 & 0 & 0 & 0 \\
0 & \beta & 0 & \frac{1}{2}\beta_y \\
0 & 0 & 0 & 0 
\end{pmatrix}\\
&\qquad
+ (t^2 - 1)
\begin{pmatrix}
0 & \frac{1}{2}V & 0 & -\frac{1}{2}(V_y - \beta\gamma_x) \\
0 & 0 & 0 & \beta\gamma \\
0 & 0 & 0 & \frac{1}{2}V \\
0 & 0 & 0 & 0 
\end{pmatrix}
\\
&\qquad
+ (t^3 - 1)
\begin{pmatrix}
0 & 0 & 0 & \frac{1}{2}\beta W \\
0 & 0 & 0 & 0 \\
0 & 0 & 0 & 0 \\
0 & 0 & 0 & 0 
\end{pmatrix}
\end{align*}
and
\begin{align*}
\frac{\partial}{\partial y}(A_t - A_1)^T &= (t - 1)
\begin{pmatrix}
0 & 0 & -\frac{1}{2}\gamma_x  & -\frac{1}{2}\gamma_{xx} \\
0 & 0 & \gamma & \frac{1}{2}\gamma_x \\
0 & 0 & 0 & 0 \\
0 & 0 & 0 & 0 
\end{pmatrix}\\
&\qquad
+ (t^2 - 1)
\begin{pmatrix}
0 & 0 & \frac{1}{2}W & \frac{1}{2}(W_x - \gamma\beta_y) \\
0 & 0 & 0 & \frac{1}{2}W \\
0 & 0 & 0 & \beta\gamma \\
0 & 0 & 0 & 0 
\end{pmatrix}
\\
&\qquad
+ (t^3 - 1)
\begin{pmatrix}
0 & 0 & 0 & \frac{1}{2}\gamma V \\
0 & 0 & 0 & 0 \\
0 & 0 & 0 & 0 \\
0 & 0 & 0 & 0 
\end{pmatrix}.
\end{align*}
We want to describe this family $d_t$ of connections invariantly. To see how to do this, let $g$ be the symmetric bilinear form defined w.r.t. to the same frame as
\begin{equation}
g:=
\begin{pmatrix}
0 & 0 & 0 & 1 \\
0 & 0 & -1 & 0 \\
0 & -1 & 0 & 0 \\
1 & 0 & 0 & \beta\gamma \\
\end{pmatrix}.
\end{equation}
This is the congruence of Lie quadrics, see e.g. Sasaki \cite{Sasaki2006}. A routine calculation (assisted by the fact that it is block upper triangular and symmetric w.r.t. $g$) shows that $\cNg$ in matrix form is
\begin{align*}
\cN^g &= 
\begin{pmatrix}
0 & -\frac{1}{2}\beta_y & 0 & \frac{1}{2}(V_y - \beta_{yy} + \beta W - 2\beta\gamma_x - \gamma\beta_x) \\
0 & 0 & 0 & 0 \\
0 & \beta & 0 & \frac{1}{2}\beta_y \\
0 & 0 & 0 & 0 
\end{pmatrix}dx\\
&\qquad
+
\begin{pmatrix}
0 & 0 & -\frac{1}{2}\gamma_x  & \frac{1}{2}(W_x-\gamma_{xx}+\gamma V - 2\gamma\beta_y - \beta\gamma_y) \\
0 & 0 & \gamma & \frac{1}{2}\gamma_x \\
0 & 0 & 0 & 0 \\
0 & 0 & 0 & 0 
\end{pmatrix}dy.
\end{align*}
We conclude that
\begin{equation}
d_t = \cD^g + t\cN^g + (t^2-1)(\chi^g + d\tau) + (t^3 - t)\psi^g,
\end{equation}
where 
\begin{align}
	\label{eq:chi}
\chi^g &= 
\frac{1}{2}
\begin{pmatrix}
0 & V & \beta\gamma & 0 \\
0 & 0 & 0 & \beta\gamma \\[0.2em]
0 & 0 & 0 & V \\
0 & 0 & 0 & 0
\end{pmatrix}dx +
\frac{1}{2}
\begin{pmatrix}
0 & \beta\gamma & W & 0 \\
0 & 0 & 0 & W  \\
0 & 0 & 0 & \beta\gamma \\
0 & 0 & 0 & 0
\end{pmatrix}dy,
\end{align}
\begin{align}
	\label{eq:psi}
\psi^g &=
\begin{pmatrix}
0 & 0 & 0 & \frac{1}{2}\beta W\\
0 & 0 & 0 & 0\\
0 & 0 & 0 & 0\\
0 & 0 & 0 & 0
\end{pmatrix}dx +
\begin{pmatrix}
0 & 0 & 0 & \frac{1}{2}\gamma V\\
0 & 0 & 0 & 0\\
0 & 0 & 0 & 0\\
0 & 0 & 0 & 0
\end{pmatrix}dy
\end{align}
and
\begin{equation}
	\label{eq:tau}
\tau = 
\begin{pmatrix}
0 & 0 & 0 & \frac{1}{2}\beta\gamma\\
0 & 0 & 0 & 0\\
0 & 0 & 0 & 0\\
0 & 0 & 0 & 0
\end{pmatrix}
\end{equation}

It will be convenient to remove the $\tau$. To do this, we gauge the family of connections using the following formula derived in \cite{Burstall2010} (see also \cite[\S8]{Burstall2004}) from the standard right logarithmic derivative formula: for $\tau\in\fsec{\tg}$, and affine connection $\nabla$ on $\tg$,
\begin{equation} \label{eq:log}
\exp(\tau)\cdot\nabla = \nabla - \nabla\tau - \frac{1}{2!}[\tau,\nabla\tau] - \frac{1}{3!}\ad^2(\tau) \nabla\tau - \dots.
\end{equation}
Now $[\chi,\tau] = 0$, so we have $\exp((t^2-1)\tau)\cdot d_t = \cD^g + t\cN^g + (t^2-1)\chi^g + (t^3-t)\psi^g$. This motivates the definition with which we begin the next section.
\section{Definition and zero-curvature formulation} \label{sec:definition}
In this section we define M\"obius-flat hypersurfaces in projective space and show that this is equivalent to the existence of a family of flat connections of a certain form. Hence we obtain a spectral deformation. 
\begin{definition}
	\label{def:mf}
Let $E$ be Legendre. If there exist one-forms $\chi^g\in\fone{\tg_{-1}\cap\h},\psi^g\in\fone{\tg_{-2}}$ for some enveloped metric $g$ on $\trR^{n+1}$ so that 
\begin{equation}\label{eq:mf_a}
\dDg\chi = \RDg
\end{equation}
 and
\begin{equation}\label{eq:mf_b}
d\psi^g + [\cNg\wedge\chi^g] = 0
\end{equation}
then we say that $(E,\chi,\psi)$ is \emph{M\"obius-flat}.
\end{definition}
If such $\chi^g$, $\psi^g$ exist for some $g$, then we can find them for any enveloped $g$ (see \S\ref{sec:gauge_theory}). This justifies viewing them as functions depending on $g$.

We saw already in \S\ref{sec:newbit} that surfaces satisfying $\beta_{yyy} = \gamma_{xxx}$ are projectively M\"obius-flat. We will see later in \S\ref{sec:conformally_flat} that the projectively M\"obius-flat hypersurfaces for $n>3$ are exactly the hypersurfaces with conformally flat second fundamental form. 

\subsection{Spectral deformation}
Now let $g$ be any enveloped metric, and suppose $\chi^g\in\fone{\tg_{-1}\cap\h}, \psi^g\in\fone{\tg_{-3}}$. Let us contemplate the family of connections $\dgt := \cDg + t\cNg + (t^2-1)\chi^g + (t^3 - t)\psi^g$. We calculate the curvature
\begin{align*}
R^{d^g_t} &= \RDg + t\dDg\cNg +  (t^2-1)\dDg\chi^g + (t^3 - t)\dDg\psi^g \\ 
&\qquad + \frac{t^2}{2}([\cNg \wedge \cNg] + (t^2-1)^2[\chi^g\wedge \chi^g] + (t^3-t)^2[\psi^g\wedge\psi^g])\\ 
&\qquad + t(t^2-1)[\cNg\wedge\chi^g] + t(t^3-t)[\cNg\wedge\psi^g] + (t^2-1)(t^3-t)[\chi^g\wedge\psi^g])\\
&=\RDg +  (t^2-1)\dDg\chi^g + (t^3 - t)(\dDg\psi^g + [\cNg\wedge\chi^g]) - \frac{t^2}{2}\RDg
\end{align*}
where we used the Gauss-Codazzi-Ricci equations $\RDg+\frac{1}{2}[\cNg\wedge\cNg] = 0$ and $d\cNg = 0$, and the fact that $[\chi^g\wedge\chi^g]$, $[\cNg\wedge\psi^g]$, $[\chi^g\wedge\psi^g]$ and $[\psi^g\wedge\psi^g]$  all vanish for algebraic reasons: the first two take values in $\tg_{-2}\cap\h = 0$ and the second two take values in $\tg_{-3} = 0$. We then observe that flatness of $\dgt$ for all $t\in\R$ is equivalent to the equations \eqref{eq:mf_a},\eqref{eq:mf_b} above.

From this we obtain a spectral deformation. In the \emph{passive} viewpoint, we replace $d$ with the flat connection $\dgt$ in the definition \ref{def:mf}; the decomposition of $d$ induced by $g$ becomes $\dgt = \cDg_t + \cNg_t$ where
\begin{align}
\cDg_t &= \cDg + (t^2-1)\chi^g\\
\cNg_t &= t\cNg + (t^3-t)\psi^g.
\end{align}
Note that $d|_{E_1} = \dgt|_{E_1}$ so that $E$ remains Legendre, and $\psi|_{E_1} = 0$ so $E$ still envelopes $g$.
Now the connection
\begin{align*}
d^g_{st} &= \cDg + st\cNg + ((st)^2-1)\chi^g + ((st)^3 - st)\psi^g\\
&= \cDg_t + s\cNg_t + (s^2 - 1)t^2\chi^g + (s^3 - s)t^3\psi^g
\end{align*}
is flat for all $s\in\R$, so $(E,t^2\chi,t^3\psi)$ is M\"obius-flat with respect to $\dgt$.

A \emph{gauge transformation} is a section $\Phi\in\fsec{\submfd\times\GL(\R^{n+1})}$; it acts on connections by $\Phi\cdot\nabla:=\Phi\circ\nabla\circ\Phi^{-1}$.
In the \emph{active} viewpoint, we find a gauge transformation $\Phi_t$ so that $\Phi_t\dgt = d$, then the spectral deformation is $\Phi_t E$.

\subsection{Gauge freedom} \label{sec:gauge_theory}
Sections $\nu$ of $\tg_{-1}\cap\m$ act freely transitively on the set of enveloped metrics by $g\mapsto\exp(-\nu)g$. Equivalently, we replace $d$ by $\exp(\nu)\cdot d$, which we may compute using the gauging formula \eqref{eq:log} from \cite{Burstall2010}. From this we see that $\cDg$, $\cNg$ become
\begin{align*}
\exp(\nu)\cdot\cDg &= \cDg - \cNg\nu - \frac{1}{2}[\nu,\cDg\nu],\\
\exp(\nu)\cdot\cNg &= \cNg - \cDg\nu - \frac{1}{2}\ad(\nu)^2(\cNg)  - \frac{1}{6}\ad(\nu)^2(\cDg\nu)
\end{align*}
respectively.
We calculate,
\begin{align*}
\exp(t\nu)\cdot\dgt &= \dgt - t\cDg\nu - t^2[\cNg,\nu] - (t^3 - t)[\chi^g,\nu] \\
&\qquad- \frac{1}{2}(t^2[\nu,\cDg\nu] - t^3\ad(\nu)^2(\cNg)) - \frac{1}{6}t^3\ad(\nu)^2(\cDg\nu)\\
&= \cDg - [\cNg,\nu] - \frac{1}{2}[\nu,\cDg\nu]\\
&\qquad + t\left(\cNg - \cDg\nu - \frac{1}{2}\ad(\nu)^2(\cNg) - \frac{1}{6}\ad(\nu)^2(\cDg\nu)\right) \\
&\qquad + (t^2 - 1)\left(\chi^g - [\cNg,\nu] - \frac{1}{2}[\nu,\cDg\nu]\right)\\
&\qquad + (t^3 - t)\left(\psi^g - \frac{1}{2}\ad(\nu)^2(\cNg) - \frac{1}{6}(\ad(\nu)^2(\cDg\nu) - [\chi^g,\nu]\right).
\end{align*}
Using the zero-curvature formulation, this justifies the claim made in the definition above. It also gives us a kind of `permutability' theorem: gauging by $t\nu$ the spectral deformation gives the spectral deformation gauged by $\nu$.

The Codazzi equation implies that the trace of $\cNg$ is closed. Since $\exp(\alpha\id)\cdot\cDg = \cDg$ and $\exp(\alpha\id)\cdot\cNg = \cNg - (d\alpha)\id$, we may assume locally, after rescaling, that $g$ is unimodular. 
\section{Examples} \label{sec:examples}

Before we look at our first example, let us recall some more theory from \cite{Burstall2004}.
A non-degenerate quadric in $\RP^n$ is a non-degenerate symmetric bilinear form on $\R^{n+1}$ up to rescaling, and so a map from $\Sigma$ in to the space of quadrics corresponds to the conformal class of a metric on $\trR^{n+1}$. Demanding that the metric be unimodular fixes it up to constant scale (using the remark at the end of \S\ref{sec:gauge_theory}). Geometrically, the condition that $E$ envelopes $g$ means that the surface and the quadric have second-order contact, in particular the asymptotic directions of the surface and the quadric at the point of contact coincide.

\subsection{Second-order envelope of a congruence of quadrics}
Akivis--Konnov \cite{Akivis1993} proved that a second order envelope of a one-parameter family of quadrics in $\RP^{n+1}$, $n\geq 4$ is conformally flat. We now extend this by showing that the envelope is M\"obius-flat for all $n \geq 3$.

So let $g$ be a metric representing a congruence of quadrics enveloping $\Lambda$, and suppose $\cNg$ has co-dimension one kernel. By the Gauss-Ricci equation we have $\RDg = 0$, so $\Lambda$ is M\"obius-flat with $\chi^g = 0 = \psi^g$.

Burstall--Calderbank \cite{Burstall2010} showed that channel surfaces in $S^3$ are M\"obius-flat, and we have just obtained the projective version of this result.

\begin{remark}
Whenever a M\"obius-flat hypersurface envelopes a curved flat, the dressing transform \cite{Terng2000} for curved flats induces a transform of enveloped hypersurfaces (though the result may not be immersed), as is easily seen using its gauge-theoretic formulation \cite{Burstall2011, Burstallb}. This serves as a partial analogue of the Eisenhart transform \cite{Eisenhart1914} of Guichard and channel surfaces, and its generalisation to conformally M\"obius-flat submanifolds (of arbitrary co-dimension) by Burstall--Calderbank \cite{Burstallb}.

In the present situation the condition of $\cNg$ having co-dimension one kernel is preserved.
\end{remark}

\subsection{Hypersurfaces with flat centro-affine metric}\label{sec:flat_ca}
Let us recall the basics of centro-affine differential geometry (e.g. from \cite{Nomizu1994, Shirokov1959, Ferapontov2004}). If we choose a linear subspace $V\subseteq\R^{n+1}$ of dimension $n$, and a vector $p\in\R^{n+1}$, then we may identify points $v\in V$ with points $\vspan{v+p}$ of an affine hyperplane in $\RP^n$. In this way we may view centro-affine geometry as a subgeometry of projective geometry. An immersion $r:\submfd \rightarrow V$ is said to be \emph{centro-affine} if $\submfd\times V = dr(T\submfd) \oplus \vspan{r}$. The trivial connection $d$ decomposes as $d_X Y = \tilde\nabla_X Y  + \hat{g}(X,Y)r$; here $\hat{g}$ is a symmetric bilinear form, the induced \emph{centro-affine metric}. The difference tensor $h:=\tilde\nabla - \hat\nabla$ where $\hat\nabla$ is the Levi-Civita connection for $\hat{g}$ is totally symmetric by the Codazzi equation, and the cubic form so defined is called the \emph{centro-affine (Fubini--Pick) cubic form}. The trace of $h$ with respect to $\hat{g}$ is called the \emph{Chebyshev covector} denoted by $T$, and \emph{proper affine hyperspheres} are characterised by the condition $T\equiv 0$.

Ferapontov showed \cite{Ferapontov2004} that if $\hat{g}$ is flat then the surface is projectively M\"obius-flat; for $n=3$ the proof consisted of giving $\beta, \gamma, V$ and $W$ in terms of a potential for the centro-affine cubic form.

To see this in our setup, we introduce a metric $g$ on $\trR^{n+1}$ adapted to the problem. First set $\Lambda = \vspan{R}$ and $\hat\Lambda = \vspan{\hat{R}}$ with $R:=p+r$ and $\hat{R}:=p-r$. Then define (as in \cite{Burstalla}) an enveloped metric using the Weyl structure $\Lambda\oplus U\oplus \hat\Lambda$ with $U:=dr(T\submfd)$. More explicitly, $\Lambda$, $\hat\Lambda$ are null, $U = (\Lambda\oplus\hat\Lambda)^\bot$ and $g(d_X R, d_Y R):=-g(d_Xd_Y R,R)$. We also declare $g(r,r) = -g(p,p) = 1$ and $g(p,r) = 0$ (N.B. such $g$ are not unimodular in general). It follows that $\cNg|_{\Lambda\oplus\hat\Lambda} = 0$, from which we see that for $D$ the induced Weyl derivative, $DR = 0 = D\hat{R}$. Now $g(\cDg R,\cDg R)$ is the centro-affine metric, so the induced Weyl connection $D^{\cDg}$ is the Levi-Civita connection $\hat\nabla$.

Now decompose $\cDg = D - \beta - \hat\beta$ using the Weyl structure, i.e. $\beta\in\fone{\stab(\hat\Lambda)^\bot}$, $\hat\beta\in\fone{\stab(\Lambda)^\bot}$. Since $d\beta = 0$ ($\cDg$ is torsion free), \eqref{eq:mf_a} holds with $\chi^g:= -\hat\beta$. Furthermore, the Gauss equation reads $d^{\cDg}\cNg - [\cNg\wedge\beta] - [\cNg\wedge\hat\beta] = 0$, the terms having weight 0, 1, -1 respectively. Thus \eqref{eq:mf_b} holds with $\psi^g = 0$.

Recall that we say a polynomial $q$ in $\fsec{\R^{n+1}}[t]$ is \emph{conserved} \cite{Burstallb} if $\dgt q(t) = 0$ for all $t\in\R$. In the present situation, $\hat{R} + t^2R$ is conserved. Indeed,
\begin{align*}
\dgt(\hat{R} + t^2 R) &=  (\cDg + t\cNg - (t^2-1)\hat\beta)(\hat{R} + t^2 R)\\
&= (\cDg + \hat\beta)\hat{R} + t\cNg\hat{R} + t^2(\cDg R - \hat\beta \hat{R}) \\
&\qquad + t^3\cNg R - t^4\hat\beta R
= 0.
\end{align*}

\begin{remark}
Let $\rho$ denote the isomorphism $\trR^{n+1} \rightarrow (\trR^{n+1})^*$ induced by $g$. Then $q\in\Gamma(\trR^{n+1})[t]$ is conserved if and only if $(\rho q)(t) := \rho(q(-t))$ is. This is because $\rho\cdot\dgt:=\rho\circ\dgt\circ\rho^{-1} = d^g_{-t}$, as connections on $(\trR^{n+1})^*$.
\end{remark}

Suppose conversely that (the contact lift of) $\Lambda$ is M\"obius-flat and there is a quadratic conserved quantity $q(t) = \hat{R} + t^2R$ with $R\in\fsec{\Lambda}$, $\hat\Lambda:=\vspan{\hat{R}}$ null and $g(R,\hat{R})>0$.
We have 
\begin{align*}
0 &= (\cDg + t\cNg + (t^2-1)\chi^g + (t^3-t)\psi^g)(\hat{R} + t^2 R)\\
&= (\cDg - \chi^g)\hat{R} + t(\cNg - \psi^g)\hat{R} + t^2(\cDg R + \chi^g \hat{R}) \\
&\qquad + t^3(\cNg R + \psi^g \hat{R}) - t^4\chi^g R + t^5\psi^g R.
\end{align*}
If we let $p$ be the constant $\frac{1}{2}(R+\hat{R})$ then the $t$ and $t^3$ terms imply $\cNg p = 0$, i.e. $g(p,\cdot)$ is constant. In particular, $g(p,p)$ and $V:=\ker(\rho q)(1)$ are constant; after rescaling we may thus assume that $g(p,p) = - 1$.

Thanks to the enveloping assumption, the condition $g(\cNg R,\hat{R}) = g(R,\cNg\hat{R}) = - g(R,\cNg R) = 0$ implies that $\cNg|_{\Lambda\oplus\hat\Lambda}$ = 0. Therefore $g$ is the adapted metric defined above, the $t$ component implies that $\psi^g = 0$, and the component constant in $t$ implies that $\chi^g = -\hat\beta$. Flatness of $\cDg - \chi^g$ thus gives flatness of $R^D$. That is to say, the centro-affine metric is flat.

Let us consider how a hypersurface with such conserved quantity transforms under the spectral deformation. The quantity $q_t(s):=q(st)$ is conserved when $d$ is replaced by $\dgt$. If $\Phi_t\cdot \dgt = d$, then $\Phi_t q_t$ is conserved for the deformed hypersurface $\Phi_t \Lambda$. However there is an arbitrary choice of constant of integration in $\Phi_t$, thus we can fix things up so that $\Phi_t p = p$ since $dp=0$. Then $\Phi_t\Lambda$ has flat metric in the \emph{same} centro-affine subgeometry.

In fact more can be said. In the passive viewpoint, if $t\neq 0$ we rescale the conserved quantity to $\hat{q}_t(s) = \frac{1}{t}\hat{R} + s^2tR$, so that $g$ is the adapted metric. Thus the centro-affine metric is scaled by $t^2$ (the change in connection may be ignored since the solder forms are the same). We already know that the cubic form is scaled by $t$. So we conclude that this spectral deformation coincides with the one given by Ferapontov \cite{Ferapontov2004} to characterise this class of hypersurfaces. He notes that for $t = 0$ one obtains a hyperplane, and this coincides with the fact in our setup that we have a constant conserved quantity.

In \S\ref{sec:pcq_surface}, we will see for the surfaces case how to formulate this in the setup of \eqref{eq:cmf_a}-\eqref{eq:cmf_c}.
\section{The main theorem} \label{sec:conformally_flat}
Having now studied projectively M\"obius-flat surfaces, it is time to address the main claim of this paper: that projectively M\"obius-flat surfaces are exactly the surfaces with flat asymptotic conformal structure when $m\geq 3$, and surfaces satisfying \eqref{eq:cmf_a}-\eqref{eq:cmf_c} when $m = 2$.

In the following subsection, we recall the tools that we shall use to analyse the equations \eqref{eq:mf_a}, \eqref{eq:mf_b}.
\subsection{Conformal Cartan Geometries and the BGG calculus}
In the case $n = 3$ the projective hypersurface inherits not only a conformal structure, but also a \emph{M\"obius structure} in the sense of \cite{Calderbank1998}. To handle this, rather than studying the conformal structure we instead use conformal Cartan geometries. In this subsection we recall the linear viewpoint, following \cite{Burstall2010}.
\begin{definition}[\cite{Bailey1994,Burstall2010}]
A \emph{conformal Cartan connection} on an $m$-manifold $\submfd$ is 
\begin{itemize}
\item{A rank $m+2$ vector bundle $V$ over $\Sigma$ equipped with a non-degenerate symmetric bilinear form $g$}
\item{An oriented null line subbundle $\Lambda \leq V$}
\item{A metric connection $\cD$ on $V$ satisfying the \emph{Cartan condition} that the algebraic map $\beta:T\Sigma\rightarrow \Hom(\Lambda,\Lambda^\bot/\Lambda), X\mapsto(\sigma\mapsto(d_X\sigma \mod \Lambda))$ is an isomorphism}.
\end{itemize}
We call this a \emph{conformal Cartan geometry} if $\cD$ is \emph{strongly torsion free} i.e. $R^\cD|_\Lambda = 0$.
\end{definition}
Restriction gives an isomorphism $\Hom(\Lambda,\Lambda^\bot/\Lambda)\cong \h/\stab(\Lambda)$, where $\h:=\so(V)$ and $\stab(\Lambda):=\{X\in\so(V)|X\Lambda\in\Lambda\}$, which yields the usual identification of $\h/\stab(\Lambda)$ with $T\submfd$ for homogeneous spaces. With respect to the Killing form, the dual of $T\submfd$ is $\stab(\Lambda)^\bot = \{X\in\so(V)|X\Lambda = 0 \text{ and } X\Lambda^\bot\subseteq\Lambda\}$. We use this without further remark in the sequel.

A conformal Cartan connection induces a conformal structure (i.e., a section of $S^2T^*\submfd \Lambda^{-2}$) $\cs$ given by $\cs(X,Y)\sigma^2 = g(\cD_X\sigma,\cD_Y\sigma)$. A \emph{Weyl structure} on $\submfd$ is a choice $\hat\Lambda$ of null line subbundle non-orthogonal to $\Lambda$. From the induced direct-sum decomposition of $V$ we obtain from $\cD$, by inclusion and projection, a connection $D$ on $\Lambda$ called the \emph{Weyl derivative}, and a connection $D^\cD$ on $T\submfd$ called the \emph{Weyl connection}; the later is conformal and torsion free so that if $D\sigma = 0$, $\sigma\in\fsec{\Lambda}$, then it is the Levi-Civita connection for the (pseudo-)Riemannian metric $\cs\sigma^2$.

As in \cite{Burstall2010} (and the more general theory set down in \cite{vCap2001}), we study conformal Cartan geometries by using concepts from Lie algebra homology, the setup of which we now recall. For $W$ a bundle that is fibrewise a representation of $T^*\submfd$, we define an operator $\del:\Omega^{k+1}(W)\rightarrow\Omega^{k}(W)$ by
\begin{equation*}
(\del\alpha)_{X_1,\ldots,X_{k-1}} := \sum_{1\leq i \leq n} \epsilon_i\cdot\alpha_{e_i,X_1,\ldots,X_{k-1}}
\end{equation*}
for any choice of dual bases $\epsilon_1,\ldots,\epsilon_n\in T^*\Sigma$ and $e_1,\ldots,e_n\in T\submfd$. This satisfies $\del^2 = 0$; the $k$-th homology group is written as $\HkW$. We define the \emph{quabla} operator $\quabla:\Omega^k(\Sigma,W)\rightarrow \Omega^k(\Sigma,W)$ by $\quabla:=d^\cD\del + \del d^\cD$; this is invertible on the image of $\del$. The kernel of $\quabla$ contains a unique representative $\Pi\alpha$ of each homology class $[\alpha]$ where $\Pi:= id - \quabla^{-1}\del d^\cD - d^\cD \quabla^{-1} \del$ satisfies $\Pi|_{\im\del} = 0$ and $\del\circ\Pi = 0$. This is Calderbank--Diemer's description \cite{Calderbank2001} of the differential lift of \v{C}ap et al. \cite{vCap2001}.

Each homology group $\HkW$ is a representation of $\stab(\Lambda)$, with $\stab(\Lambda)^\bot$ acting trivially. There is a unique element in $\stab(\Lambda)/\stab(\Lambda)^\bot$, called the grading element, that acts as $-1$ on $\stab(\Lambda)^\bot$, $0$ on $\stab(\Lambda)/\stab(\Lambda)^\bot$ and $1$ on $\h/\stab(\Lambda)$ (c.f. \cite{vCap2009,Calderbank2005}). When this element acts as the scalar $i$ on $\HkW$, we say that $\HkW$ has \emph{weight} $i$.

Calderbank--Diemer \cite{Calderbank2001} defined multi-linear differential operators between homology bundles which we shall exploit. Given a bilinear pairing $\wedge:W_1\times W_2 \rightarrow W_3$ of representations, we define $\bggCup:\fsec{H^k(\submfd,W_1)}\times \fsec{H^l(\submfd,W_1)} \rightarrow \fsec{H^{k+l}(\submfd,W_3)}$ by $[\omega]\bggCup[\tau] := [\Pi((\Pi\omega)\wedge(\Pi\tau))]$.

When $\cD$ is strongly torsion free, there is a unique section $Q\in\fsec{S^2 T^*M}\cap\im\del$ such that $\cD - Q$ is \emph{normal} (i.e. $\del R^{\cD - Q} = 0$). On a surface this means that $Q$ is trace-like. A conformal Cartan geometry is determined entirely (up to isomorphism) by the conformal structure when $m \geq 3$. For $m = 2$, one also needs to know an additional piece of structure induced by the connection, called the \emph{M\"obius structure} \cite{Calderbank1998, Burstall2010}.

The Weyl curvature of $\cs$ is $W=\pi R^{\cD - Q}$ where $\pi:\stab(\Lambda)\rightarrow\stab(\Lambda)/\stab(\Lambda)^\bot$ is the projection; of course when $m = 2,3$ this vanishes automatically for algebraic reasons. When $m = 2,3$, the Cotton-York curvature $\CYc$ is $R^{\cD - Q}$. This is the extension to the case $m=2$ provided by \cite{Calderbank1998}.

Recall \cite{Burstall2004} that, a unimodular metric $g$ enveloped by $E$ is said to be in the \emph{Darboux family} when $\cNg|_{E_1} = 0$. Given $g$ in the Darboux family, we obtain a conformal Cartan geometry, namely $(\trR^{n+1},g,E_1,\cDg)$. There is a unique unimodular metric $g$ satisfying the condition $\del\cNg = 0$, which  we call \emph{normal}; in fact this is the Lie quadric congruence \cite{Burstall2004}. This is the conformal Cartan geometry  that we shall use.

\subsection{Spectral deformation}
As a simple application of this, let us check that the spectral deformation scales the Darboux cubic form.
Suppose $E$ is M\"obius-flat. We work in the passive viewpoint. The solder form of $\cDg_t$ is equal to that of $\cDg$, so the conformal structure remains unchanged. Since $\del\cNg_t = t\del\cNg + (t^3 - t)\del\psi^g = 0$, the normal unimodular metric is still $g$. So $\dcf$, the homology class (see \cite{Burstall2004}) of $\cNg$ becomes $[\cNg_t] = t[\cNg] =  t\dcf$.

\subsection{The equivalence, part one}
Equations \eqref{eq:mf_a}, \eqref{eq:mf_b} can be reduced to a form involving the Darboux cubic form, a quadratic differential, the Cotton-York curvature, and one of the first order differential parings of Calderbank-Diemer.

To effect this, we use an analysis very similar to that in the proof of \cite[Theorem 16.3]{Burstall2010}.

\begin{proposition}
	\label{prop:mf_homo}
	For $m\geq 3$, $E$ is M\"obius-flat if and only if it has conformally flat second fundamental form. For $m=2$, $E$ is M\"obius-flat if and only if there exists $q\in\fsec{S^2_0 T^*\Sigma}$ with
\begin{equation}
	\label{eq:mf_homo_a}
	dq = \CYc
\end{equation}
and
\begin{equation}
	\label{eq:mf_homo_b}
	q\bggCup\dcf = 0. 
\end{equation}
\end{proposition}
Here $dq:=d^Dq$, independently of choice of Weyl connection.
\begin{proof}
Let $g$ be the normal unimodular metric. Set $q:=\chi^g - Q$, where $\cD = \cDg - Q$ is the normal connection. If $\cDg-\chi^g = \cD - q$ is flat, then $\del d^{\cD}q = 0$, thus $q$ is the differential lift of its homology class. Now when $m = 2$, the homology $H_1(T^*\submfd,\h)$ for $m = 2$ has weight $-2$ and consists of symmetric trace-free bilinear forms so $q$ is a quadratic differential. When $m \geq 3$, the homology has weight $0$, so $q$ vanishes. Thus \eqref{eq:mf_a} is equivalent to $dq = \CYc$ when $m = 2$, and $0 = R^{\cD}$ when $m \geq 3$. Now for $m = 3$ this is the same as $C = 0$; for $m > 3$ this is equivalent to the vanishing of the Weyl curvature since $R^{\cD} = \Pi W$ (using the Bianchi identity, and that $\cD$ is normal) \cite{Burstall2004}.

Since $\psi^g$ is in the image of $\del$, we may rewrite \eqref{eq:mf_b} as $0 = (id - \dDg\quabla^{-1}\del)[\cNg\wedge\chi^g] = 0$. Using the Codazzi equation and the Jacobi identity, $\dDg[\cNg\wedge Q] = 0$; also when $m = 2$, there are no 3-forms, so similarly $\dDg[\cNg\wedge q] = 0$. Therefore we may rewrite \eqref{eq:mf_b} as $\Pi[\cNg\wedge q] = 0$ when $m = 2$, and as $\Pi[\cNg\wedge Q] = 0$ when $m = 3$. The latter equation is in fact automatic since the homology $H_2(T^*\submfd,\m)$ has weight $-2$ in that case \cite{Burstalla}. 

To conclude, if $(E,\chi,\psi)$ is M\"obius flat, then when $m = 2$, $q$ satisfies the required conditions, and when $m = 3$ we have that $(\submfd,\cs)$ is conformally-flat. Conversely, if $m = 2$ and $q$ satisfies $dq = \CYc$, $q\bggCup\dcf = 0$ then $E$ is M\"obius-flat with $\chi^g = Q + q$, otherwise if $m > 2$ and $(\submfd,\cs)$ is conformally-flat then $E$ is M\"obius-flat with $\chi^g = Q$. In either case, $\psi^g = \quabla^{-1}\del[\cNg\wedge\chi^g]$.
\end{proof}

This condition on $q$ was used previously by Burstall--Calderbank to give a similar reformulation of projective and Lie applicability (private communication).
\subsection{The equivalence, part two}
\label{sec:hom_surfaces}

We are now going to see that, for the case of a surface, equations \eqref{eq:mf_a} and \eqref{eq:mf_b} are in fact co-ordinate invariant versions of \eqref{eq:cmf_a} and \eqref{eq:cmf_b},\eqref{eq:cmf_c} respectively.
To assist in these calculations, we let $e_1,e_2$ be the co-ordinate vector fields for $x,y$ respectively, with dual basis $\epsilon_1,\epsilon_2$. Also let $\cs$ be conformal structure in induced by the conformal Cartan geometry $(\tR^{n+1}, g, E_1, \cDg)$. We work with the frame $(\sigma, \sigma_x, \sigma_y, \hat\sigma)$ using $\sigma$ as in \S\ref{sec:intro} and $\hat\sigma = \sigma_{xy} - \frac{1}{2}\beta\gamma$.

It follows that 

\begin{equation*}
\epsilon_1 =
\begin{pmatrix}
0 & 1 & 0 & 0\\
0 & 0 & 0 & 0\\
0 & 0 & 0 & 1\\
0 & 0 & 0 & 0
\end{pmatrix},
\quad
\epsilon_2 = 
\begin{pmatrix}
0 & 0 & 1 & 0\\
0 & 0 & 0 & 1\\
0 & 0 & 0 & 0\\
0 & 0 & 0 & 0
\end{pmatrix}
\end{equation*}
and that $e_1, e_2$ are parallel for $D$.

First let us deal with \eqref{eq:mf_a}. For the right hand side, we need the curvature of $\cD$ which in turn comes from $R^{\cDg}$ and $Q$. 
\begin{equation*}
R^{\cDg}_{e_1,e_2} = -\frac{1}{2}[\cNg\wedge\cNg]_{e_1,e_2} = -[\cNg_{e_1},\cNg_{e_2}]
=
\begin{pmatrix}
0 & -\frac{1}{2}\beta\gamma_x & \frac{1}{2}\gamma\beta_y & 0 \\
0 & \beta\gamma & 0 & \frac{1}{2}\beta_y\gamma \\
0 & 0 & -\beta\gamma & -\frac{1}{2}\beta\gamma_x \\
0 & 0 & 0 & 0.
\end{pmatrix}
\end{equation*}
Thus if $Q = \alpha\cs$, then
\begin{align*}
(d^{\cDg}Q)_{e_1,e_2}
&= \cD^g_{e_1}Q_{e_2} - \cD^g_{e_2}Q_{e_1}\\
&= \cD^g_{e_1}\alpha\epsilon_1 - \cD^g_{e_2}\alpha\epsilon_2\\
&= \alpha
\begin{pmatrix}
-1 & 0 & 0 & 0 \\
0 & 1 & 0 & 0 \\
0 & 0 & -1 & 0 \\
0 & 0 & 0 & 1
\end{pmatrix}
+ \alpha_x\epsilon_1 - \alpha_y\epsilon_2 - \alpha
\begin{pmatrix}
-1 & 0 & 0 & 0 \\
0 & -1 & 0 & 0 \\
0 & 0 & 1 & 0 \\
0 & 0 & 0 & 1
\end{pmatrix}\\
&=
\begin{pmatrix}
0 & \alpha_x & -\alpha_y  & 0 \\
0 & 2\alpha & 0 & -\alpha_y \\
0 & 0 & -2\alpha & \alpha_x \\
0 & 0 & 0 & 0
\end{pmatrix}\\
\end{align*}
also $[Q\wedge Q] = 0$, so $\cDg-Q$ is normal when $\alpha = \frac{1}{2}\beta\gamma$.

Now we take an arbitrary trace-free quadratic differential
\begin{equation*}
q = a\epsilon_1\otimes\epsilon_1 + b\epsilon_2\otimes\epsilon_2,
\end{equation*}
and compute the exterior derivative (coupled with the Weyl connection):
\begin{equation*}
(d^Dq)_{e_1,e_2} = D_{e_1}q(e_2) - D_{e_2}q(e_1)
= D_{e_1}(b\epsilon_2) - D_{e_2}(a\epsilon_1)
= b_x\epsilon_2 - a_y\epsilon_1.
\end{equation*}

So
\begin{align*}
(dq - R^{\cDg} - Q))_{e_1,e_2} 
&= b_x\epsilon_2 - a_y\epsilon_1 - ((-\frac{1}{2}\beta\gamma_x - \frac{1}{2}\beta_x\gamma - \frac{1}{2}\beta\gamma_x)\epsilon_1 + (-\frac{1}{2}\gamma\beta_y - \frac{1}{2}\beta_y\gamma + \frac{1}{2}\beta\gamma_y)\epsilon_2)\\
&= (\frac{1}{2}\beta_x\gamma + \beta\gamma_x - a_y)\epsilon_1
+ (-\frac{1}{2}\beta\gamma_y - \beta_y\gamma + b_x)\epsilon_2
\end{align*}
and equation \eqref{eq:mf_a} becomes 
\begin{align*}
2a_y &= 2\beta\gamma_x + \beta_x\gamma\\
2b_x &= 2\beta_y\gamma + \beta\gamma_y.\\
\end{align*}
Second we attend to \eqref{eq:mf_b}. The one-forms $\cNg$ and $q$ are the differential lifts of their homology classes, so to compute the cup product we need to calculate $(\id - d^{\cDg}\quabla^{-1}\del)[\cNg\wedge q]$ (c.f. the proof of Proposition \ref{prop:mf_homo}).
So,
\begin{align}
	[\cNg\wedge q]_{e_1,e_2} &= [\cNg_{e_1},q_{e_2}] - [\cNg_{e_2}, q_{e_1}]\notag\\
&=
\begin{pmatrix}
0 & 0 & 0 & -\frac{1}{2}\beta_y b\\
0 & 0 & 0 & 0\\
0 & 0 & 0 & \beta b\\
0 & 0 & 0 & 0
\end{pmatrix}
-
\begin{pmatrix}
0 & \beta b & 0 & \frac{1}{2}\beta_y b\\
0 & 0 & 0 & 0\\
0 & 0 & 0 & 0\\
0 & 0 & 0 & 0
\end{pmatrix}\notag\\
&\qquad-
\begin{pmatrix}
0 & 0 & 0 & -\frac{1}{2}\gamma_x a\\
0 & 0 & 0 & \gamma a\\
0 & 0 & 0 & 0\\
0 & 0 & 0 &0 
\end{pmatrix}
+
\begin{pmatrix}
0 & \gamma a & 0 & \frac{1}{2}\gamma_x a\\
0 & 0 & 0 & 0\\
0 & 0 & 0 & 0\\
0 & 0 & 0 & 0
\end{pmatrix}\notag\\
&=
\begin{pmatrix}
0 & -\beta b & \gamma a & -\beta_y b + \gamma_x a\\
0 & 0 & 0 & -\gamma a\\
0 & 0 & 0 & \beta b\\
0 & 0 & 0 & 0
\end{pmatrix}.\label{eq:surface_a}
\end{align}
Thus
\begin{align*}
(\del[\cNg\wedge q])_{e_1,e_2} = 
\begin{pmatrix}
0 & 0 & 0 & -2\beta b \epsilon_1 - 2\gamma a \epsilon_2\\
0 & 0 & 0 & 0\\
0 & 0 & 0 & 0\\
0 & 0 & 0 & 0
\end{pmatrix}
\end{align*}
and it follows that
\begin{equation}
	\label{eq:surface_b}
\psi^g:=-\quabla^{-1}\del[\cNg\wedge q] = 
\begin{pmatrix}
0 & 0 & 0 & \beta b\epsilon_1 + \gamma a \epsilon_2\\
0 & 0 & 0 & 0\\
0 & 0 & 0 & 0\\
0 & 0 & 0 & 0
\end{pmatrix}.
\end{equation}
To verify this last, we confirm that $\quabla$ applied to the right hand side is equal to $-\del[\cNg\wedge q]$.
We calculate
\begin{align*}
(d^{\cDg}\begin{pmatrix}
0 & 0 & 0 & \beta b\epsilon_1 + \gamma a \epsilon_2\\
0 & 0 & 0 & 0\\
0 & 0 & 0 & 0\\
0 & 0 & 0 & 0
\end{pmatrix})_{e_1,e_2} 
&= \cDg_{e_1} 
\begin{pmatrix}
0 & 0 & 0 & \gamma a\\
0 & 0 & 0 & 0\\
0 & 0 & 0 & 0\\
0 & 0 & 0 & 0
\end{pmatrix}
- \cDg_{e_2} 
\begin{pmatrix}
0 & 0 & 0 & \beta b\\
0 & 0 & 0 & 0\\
0 & 0 & 0 & 0\\
0 & 0 & 0 & 0
\end{pmatrix}\\
&=
\begin{pmatrix}
0 & \beta b & -\gamma a & 0\\
0 & 0 & 0 & \gamma a\\
0 & 0 & 0 & \beta b\\
0 & 0 & 0 & 0
\end{pmatrix}
\mod\wedge^2 E
\end{align*}
and so
\begin{align*}
(\del d^{\cDg} \begin{pmatrix}
0 & 0 & 0 & \beta b\epsilon_1 + \gamma a \epsilon_2\\
0 & 0 & 0 & 0\\
0 & 0 & 0 & 0\\
0 & 0 & 0 & 0
\end{pmatrix})_{e_1}
&= \epsilon_2\cdot(d^{\cDg}\begin{pmatrix}
0 & 0 & 0 & \beta b\epsilon_1 + \gamma a \epsilon_2\\
0 & 0 & 0 & 0\\
0 & 0 & 0 & 0\\
0 & 0 & 0 & 0
\end{pmatrix})_{e_2,e_1}\\
&= - \left[
\begin{pmatrix}
0 & 0 & 1 & 0\\
0 & 0 & 0 & 1\\
0 & 0 & 0 & 0\\
0 & 0 & 0 & 0
\end{pmatrix},
\begin{pmatrix}
0 & \beta b & -\gamma a & 0\\
0 & 0 & 0 & \gamma a\\
0 & 0 & 0 & -\beta b\\
0 & 0 & 0 & 0
\end{pmatrix}\right]\\
&=
\begin{pmatrix}
0 & 0 & 0 & -2\beta b\\
0 & 0 & 0 & 0\\
0 & 0 & 0 & 0\\
0 & 0 & 0 & 0
\end{pmatrix}.
\end{align*}
The calculation for the $\epsilon_2$ component is similar.

Thus \eqref{eq:mf_homo_b} becomes
\begin{align*}
(d^{\cDg}\psi + [\cNg\wedge q])_{e_1,e_2}
&=
\begin{pmatrix}
0 & \beta b & -\gamma a & (\gamma a)_x - (\beta b)_y\\
0 & 0 & 0 & \gamma a\\
0 & 0 & 0 & -\beta b\\
0 & 0 & 0 & 0
\end{pmatrix}
+
\begin{pmatrix}
0 & \beta b & \gamma a & \gamma_x a - \beta_y b\\
0 & 0 & 0 & -\gamma a\\
0 & 0 & 0 & \beta b\\
0 & 0 & 0 & 0
\end{pmatrix}\\
&=
\begin{pmatrix}
0 & 0 & 0 & -2\beta_y b - \beta b_y + 2\gamma_x a + \gamma a_x\\
0 & 0 & 0 & 0\\
0 & 0 & 0 & 0\\
0 & 0 & 0 & 0
\end{pmatrix}
\end{align*}
and hence \eqref{eq:mf_b} reads
\begin{equation*}
2\beta_y b - \beta b_y = 2\gamma_x a + \gamma a_x.
\end{equation*}

As a by product, we have proved that if equations \eqref{eq:mf_a}, \eqref{eq:mf_b} hold then 
\begin{equation*}
\chi^g = 
\begin{pmatrix}
0 & adx + \frac{1}{2}\beta\gamma dy  & \frac{1}{2}\beta\gamma dx + bdy & 0\\
0 & 0 & 0 & \frac{1}{2}\beta\gamma dx + bdy\\
0 & 0 & 0 & adx + \frac{1}{2}\beta\gamma dy\\
0 & 0 & 0 & 0
\end{pmatrix}
\end{equation*}
and
\begin{equation*}
\psi^g = 
\begin{pmatrix}
0 & 0 & 0 & \beta b dx + \gamma a dy\\
0 & 0 & 0 & 0\\
0 & 0 & 0 & 0\\
0 & 0 & 0 & 0
\end{pmatrix}.
\end{equation*}

\section{Polynomial conserved quantities for surfaces} \label{sec:pcq_surface}
Now we want to interpret the result of \S\ref{sec:flat_ca} about conserved quantities in terms of the classical set-up described in the introduction.

Throughout this section, all matrices are written with respect to the frame  $(R,R_x,R_y,\hR)$ for a fixed choice of asymptotic co-ordinates $x,y$ and then $\sigma$ is as in \S\ref{sec:intro}.
Let $\clq$ be the congruence of Lie quadrics. Let $\nu\in\fsec{\tg_{-1}}$ be such that $\exp(-\nu)(R, R_x, R_y,\hR) = e^{\alpha}(\sigma, \sigma_x, \sigma_y, -2\hat\sigma)$ where $\alpha$ is a potential for the Chebyshev covector so that $\clq = \exp(-(\nu+\alpha\id))g$; write
\begin{equation*}
\nu=
\begin{pmatrix}
0 & \frac{1}{2}a & \frac{1}{2}b & c\\
0 & 0 & 0 & b \\
0 & 0 & 0 & a \\
0 & 0 & 0 & 0
\end{pmatrix}.
\end{equation*}
We know, from \S\ref{sec:gauge_theory}, that $\exp((t-1)(\nu + \alpha\id))\cdot d^g_t$ is of the form
\begin{equation*}
d_t = \cDg + t\cNg + (t^2-1)\chi^g + (t^3-t)\psi^g
\end{equation*}
and has conserved quantity
\begin{align*}
\exp((t-1)(\nu+\alpha\id))(\hR + t^2R) 
&= e^{(t-1)\alpha}\exp(-\nu)(\hR + t^2R + t(cR + bR_x + aR_y) \\
&\qquad+ \frac{t^2}{2}(\frac{1}{2}abR + \frac{1}{2}abR))\\
&=e^{t\alpha}(-2\hat\sigma + t(c\sigma + b\sigma_x + a\sigma_y) + t^2(1 + \frac{1}{2}ab)\sigma).
\end{align*}

We now gauge both the pencil of connections and the conserved quantity by $\exp(-t\alpha\id)$ to get:
\begin{align*}
\exp(-t\alpha\id)\cdot d_t &= d_t - d_t(-t\alpha\id)\\
&= d_t + td\alpha\id\\
& = \cDg + t(\cNg + d\alpha) + (t^2-1)\chi^g + (t^3-t)\psi^g
\end{align*}
and
\begin{equation*}
-2\hat\sigma + t(c\sigma + b\sigma_x + a\sigma_y) + t^2(1 + \frac{1}{2}ab)\sigma.
\end{equation*}

Notice $(\hat\cD - \hat\chi)\hR = 0$ so $(\cD^g - \chi^g)\hat\sigma = 0$. It follows that $(\cD^g - (\chi^g + d\tau))\sigma_{xy} = 0$ where $\tau$ is given by \eqref{eq:tau} as above, and thus $\chi^g$,$\psi^g$ are given by \eqref{eq:chi},\eqref{eq:psi}. Now we unpack the conservation equation. First we work slightly more generally with $v_0\in\fsec{\tR^{n+1}}$, $v_1\in\fsec{\Lambda^\perp}$ and $v_2\in\fsec{\Lambda}$.
Separating powers of $t$, we get
\begin{equation*}
0 = (\cDg + t(\cNg + d\alpha) + (t^2-1)\chi^g + (t^3-t)\psi^g)(v_0 + tv_1 + t^2 v_2)
\end{equation*}
iff
\begin{align*}
0 &= \cDg v_0 - \chi v_0\\
0 &= \cDg v_1 + \cNg v_0 + d\alpha v_0 - \chi v_1 - \psi v_0\\
0 &= \cDg v_2 + \cNg v_1 + d\alpha v_1 + \chi v_0 - \chi v_2 - \psi v_1 \\
0 &= \cNg v_2 + d\alpha v_2 + \chi v_1 + \psi v_0 - \psi v_2\\
0 &= \chi v_2 + \psi v_1\\
0 &= \psi v_2
\end{align*}
which is, taking into account kernels, equivalent to 
\begin{align*}
0 &= \cDg v_1 + \cNg v_0 + d\alpha v_0 - \chi v_1 - \psi v_0\\
0 &= \cDg v_2 + \cNg v_1 + d\alpha v_1 + \chi v_0\\
0 &= d\alpha v_2 + \chi v_1 + \psi v_0 - \psi v_2\\
\end{align*}

Now we read off our case, with $v_0 = -2\hat\sigma$, $v_1 = c\sigma + b\sigma_x + a \sigma_y$ and $v_2 = (1 + \frac{1}{2}ab)\sigma$. Firstly in the $x$ direction:

\begin{align*}
0 &= c_x \sigma + c\sigma_x + b_x\sigma_x + \frac{1}{2}Vb\sigma + a_x\sigma_y + a\sigma_{xy} - 2(\frac{1}{2}(\beta W - \beta_{yy})\sigma + \frac{1}{2}\beta_y\sigma_y) \\
&\qquad- 2\alpha_x\hat\sigma - \frac{1}{2}Vb\sigma - \frac{1}{2}\beta\gamma a\sigma - (-2)\frac{1}{2}\beta W\sigma\\
&= (c_x + \frac{1}{2}bV - \beta W + \beta_{yy} + \alpha_x\beta\gamma - \frac{1}{2}bV - \frac{1}{2}\beta\gamma a + \beta W)\sigma\\
&\qquad+ (c + b_x) \sigma_x + (a_x - \beta_y)\sigma_y + (a - 2\alpha_x) \sigma_{xy}\\
0 &= \sigma_x + b(-\frac{1}{2}\beta_y\sigma + \beta\sigma_y) + \alpha_x c\sigma + \alpha_x b\sigma_x + \alpha_x a\sigma_y - 2(\frac{1}{2}\beta\gamma\sigma_x + \frac{1}{2} V\sigma_y)\\
&= (-\frac{1}{2}\beta_yb + \alpha_x c)\sigma + (1+ \alpha_xb  -\beta\gamma)\sigma_x + (b\beta + \alpha_x a - V)\sigma_y\\
0 &= \alpha_x\sigma + (\frac{1}{2}Vb + \frac{1}{2}\beta\gamma a)\sigma + (-2)\frac{1}{2}\beta W\sigma\\
&= (\alpha_x + \frac{1}{2}(Vb + \beta\gamma a) - \beta W)\sigma.
\end{align*}

Secondly in the $y$-direction:
\begin{align*}
0 &= c_y\sigma + c\sigma_y + b_y\sigma_x + b\sigma_{xy} + a_y \sigma_y + \frac{1}{2}Wa\sigma + (-2)(\frac{1}{2}(\gamma V - \gamma_{xx})\sigma + \frac{1}{2}\gamma_x \sigma_x)\\
&\qquad - 2\alpha_y\hat\sigma - (\frac{1}{2}\beta\gamma b + \frac{1}{2} W a)\sigma - (-2)\frac{1}{2}\gamma V \sigma\\
&= (c_y + \frac{1}{2}W a - \gamma V + \gamma_{xx}  + \alpha_y \beta\gamma - \frac{1}{2}\beta\gamma b - \frac{1}{2} W a + \gamma V) \sigma\\
&\qquad+ (b_y - \gamma_x) \sigma_x + (c + a_y)\sigma_y + (b - 2\alpha_y)\sigma_{xy}\\
0 &= \sigma_y + a(-\frac{1}{2}\gamma_x\sigma + \gamma\sigma_x) + \alpha_y c\sigma + \alpha_y b\sigma_x + \alpha_y a \sigma_y + (-2)(\frac{1}{2}W\sigma_x + \frac{1}{2}\beta\gamma\sigma_y)\\
&= (-\frac{1}{2}\gamma_x a + \alpha_yc)\sigma + (a\gamma + \alpha_yb - W)\sigma_x + (1 + \alpha_y a - \beta\gamma)\sigma_y\\
0 &= \alpha_y \sigma + (\frac{1}{2}\beta\gamma b + \frac{1}{2}W a )\sigma + (-2)\frac{1}{2}\gamma V\sigma\\
&= (\alpha_y + \frac{1}{2}(\beta\gamma b + Wa) - \gamma V)\sigma.
\end{align*}

This is then equivalent to:
\begin{align*}
0 &= c_x + \beta_{yy} & 0&=c_y + \gamma_{xx}\\
a_x &= \beta_y & b_y &= \gamma_x\\
a &= 2\alpha_x & b &= 2\alpha_y \\
c&=-a_y & c &= -b_x\\
\alpha_x c &= \frac{1}{2}\beta_y b & \alpha_yc &= \frac{1}{2}\gamma_x a\\
0 &= 1 + \alpha_xb  -\beta\gamma& 0 &= 1 + \alpha_y a - \beta\gamma \\
V &= \beta b + \alpha_x a & W &= a\gamma + \alpha_y b\\
0 &= \alpha_x + \frac{1}{2}(Vb + \beta\gamma a) - \beta W & 0 &= \alpha_y +\frac{1}{2}(\beta\gamma b + Wa) - \gamma V
\end{align*}
which is equivalent to $a = 2\alpha_x, b = 2\alpha_y, c = -2\alpha_{xy}$ together with the equations given in Theorem \ref{thm:pcq_classical}
\section*{Funding}
This work was partly supported by an Engineering and Physical Sciences Research Council DTA.

\section*{Acknowledgements}
I would like to thank my PhD supervisor David Calderbank for suggesting this topic; Jenya Ferapontov and Udo Hertrich-Jeromin for discussion. 

\bibliographystyle{abbrvurl}
\bibliography{mf}
\end{document}